\documentclass[12pt]{article}
\usepackage[utf8]{inputenc}
\usepackage[english]{babel}
\usepackage{graphicx}
\usepackage{amssymb}
\usepackage{amsthm}
\usepackage{amsmath}

\usepackage{subcaption}
\usepackage[font=footnotesize]{caption}

\usepackage{tikz}\usetikzlibrary{matrix}
\usepackage[left=1.5in,right=1in]{geometry}

\usepackage[T1]{fontenc}

\usepackage[shadow, colorinlistoftodos,textsize=tiny]{todonotes}
\setuptodonotes{fancyline, backgroundcolor=gray!10,bordercolor=gray}

\usepackage{authblk}
\usepackage{hyperref}
\hypersetup{
	hidelinks,
  colorlinks  = true,    urlcolor    = blue,    linkcolor   = blue,    citecolor   = blue,      pdfauthor   = {Elías Mochán},pdfsubject  = {Maniplexes},pdftitle    = {Polytopality of 2-orbit maniplexes},  }

\newtheorem{teo}{Theorem}[section]

\newtheorem{coro}[teo]{Corollary}

\newtheorem{lema}[teo]{Lemma}
\theoremstyle{definition}
\newtheorem{Def}[teo]{Definition}

\newtheorem{af}[teo]{Claim}
\newenvironment{dem}[1][Proof]{\begin{trivlist}
  \item[\hskip \labelsep {\bfseries #1}]}{\end{trivlist}}

\theoremstyle{remark}

\newcommand{\Z}{\mathbb{Z}}

\newcommand{\F}{\mathcal{F}}

\newcommand{\p}{\mathcal{P}}

\newcommand{\G}{\mathcal{G}}
\newcommand{\M}{\mathcal{M}}
\newcommand{\X}{\mathcal{X}}
\newcommand{\T}{\mathcal{T}}

\newcommand{\rk}{{\bf rank}}
\newcommand{\ol}[1]{\overline{#1}}\newcommand{\os}[1]{\widetilde{#1}}

\newcommand{\gen}[1]{\langle #1 \rangle} \newcommand{\Aut}{\Gamma}
\newcommand{\supp}{{\rm \bf supp}}
\newcommand{\Fac}{{\rm Fac}}
\newcommand{\Mon}{{\rm Mon}}

\title{Polytopality of 2-orbit maniplexes
}
\author{Elías Mochán}
\affil{Department of Mathematics, Northeastern University, 02115
  Boston, USA\thanks{\tt email:j.mochanquesnel@northeastern.edu}}
\begin{document}
\maketitle

\abstract{\emph{Abstract polytopes} are a combinatorial generalization
  of convex and skeletal polytopes. Counting how many flag orbits a
  polytope has under its automorphism group is a way of measuring how
  symmetric it is. Polytopes with one flag orbit are called
  \emph{regular} and are very well known. Polytopes with two flag
  orbits (called \emph{2-orbit polytopes}) are, however, way more
  elusive. There are $2^n-1$ possible classes of 2-orbit polytopes in
  rank (dimension) $n$, but for most of those classes, determining
  whether or not they are empty is still an open problem. In 2019, in
  their article \emph{An existence result on two-orbit maniplexes},
  Pellicer, Poto\v{c}nik and Toledo constructed 2-orbit
  \emph{maniplexes} (objects that generalize abstract polytopes and
  maps) in all these classes, but the question of whether or not they
  are also polytopes remained open. In this paper we use the results
  of a previous paper by the author and Hubard to show that some of
  these 2-orbit maniplexes are, in fact, polytopes. In particular we
  prove that there are 2-orbit polytopes in all the classes where
  exactly two kinds of reflections are forbidden. We use this to show
  that there are at least $n^2-n+1$ classes of 2-orbit polytopes of
  rank $n$ that are not empty. We also show that the maniplexes
  constructed with this method in the remaining classes satisfy all
  but (possibly) one of the properties necessary to be polytopes,
  therefore we get closer to proving that there are 2-orbit polytopes
  in all the classes.}
\section{Introduction}\label{s:Intro}

Abstract polytopes are posets with some properties that generalize
those of the face-poset of convex and skeletal polytopes. Since
their introduction by Schulte in~\cite{EgonPHD}, the main approach to
study abstract polytopes has been by their symmetries. It is known
that the group of symmetries (or automorphism group) of an abstract
polytope acts freely on its \emph{flags}: chains consisting of one
element of each rank (dimension). Therefore, we may measure how
symmetric a polytope is by counting how many flag orbits it has under
the action of its automorphism group. A \emph{$k$-orbit polytope} is
one that has $k$ flag orbits under this action. The smaller the number
of flag orbits, the more symmetric the polytope. We call 1-orbit
polytopes \emph{regular}, and they have been by far the most
studied. The book~\cite{ARP} is the standard reference and it is
dedicated solely to abstract regular polytopes.

It is quite clear that 2-orbit polytopes are the second most symmetric
kind of polytopes after the regular ones, so, naturally, they have been
the second most studied family. Nevertheless, the general theory of
2-orbit polytopes has been much more challenging than that of regular
ones.

One can classify 2-orbit polytopes of rank $n$ in $2^n-1$ different
\emph{symmetry types}. These symmetry types are determined by a proper
subset $I\subset \{0,1,\ldots,n-1\}$, that in some way consists of the
``permitted kinds of reflections'', and are denoted by $2^n_I$. One of those symmetry types consists of the so called {\em chiral}
polytopes: those with all possible ``rotational'' symmetry but no
``reflection'' symmetry at all ($2^n_\emptyset$). Historically,
chiral polytopes have been the most studied type of 2-orbit
polytopes. The main theory of (abstract) chiral polytopes was
developed in 1991 by Schulte and Weiss~\cite{QuiralesEgonAsia} but the
existence of chiral polytopes in any rank was one of the main
questions to consider. Rank 3 chiral polytopes had been studied in the
context of maps on surfaces, and in the 1970's Coxeter gave examples
of rank 4 chiral polytopes arising as quotients of hyperbolic
tessellations~\cite{coxeter_twisted}. However, it is not until 2008
that algebraic methods were developed to find (finite) chiral
polytopes of rank 5~\cite{QuiralesConstruct}. In 2010, almost 20 years
after the publication of~\cite{QuiralesEgonAsia}, Pellicer proved that
there exist abstract chiral polytopes in any rank at greater than 2~\cite{Quirales}.

Also in 2010, Hubard~\cite{2OrbPol1} described a way to construct
2-orbit polyhedra of any symmetry type from groups;
examples of 2-orbit polyhedra for each of the 7 possible symmetry
types are well-known. In 2016, Pellicer~\cite{Cuboctaedro}
found geometric examples of 2-orbit polytopes of any rank with a fixed
symmetry type. However, the challenge of finding 2-orbit polytopes of
any possible type has not had that much of an advancement.

\emph{Maniplexes} are a generalization of polytopes where some
conditions on connectivity are relaxed.  In 2019, Pellicer,
Poto\v{c}nik and Toledo~\cite{2OrbMani} found a way to construct
2-orbit maniplexes of any type. To do this, they used a \emph{voltage
  assignment} on the symmetry type graph. However, the question of
whether or not their examples are flag-graphs of
polytopes remained open. This is due to the fact that there was no
method to find the intersection properties that the so called voltage
group should satisfy for the derived maniplex to be
polytopal. However, in~\cite{IntPropYo}, Hubard and the author of this
paper described a method to find the intersection properties for every
symmetry type for any number of orbits (see Theorem~\ref{t:IntProp}).
Using this method for 2-orbit maniplexes, we have a way to test and
determine their polytopality.

The 2-orbit $n$-maniplexes constructed in~\cite{2OrbMani} depend on
some choices: a regular $(n-1)$-maniplex and a \emph{monodromy}
(to be defined in~\ref{s:basic}) of that maniplex, both satisfying
some conditions.
In their paper, Pellicer, Poto\v{c}nik and Toledo proved that given a
2-vertex symmetry type, there exists a maniplex $\M$ and a monodromy
$\eta$ of $\M$ satisfying such conditions and thus, the $n$-maniplex
obtained has the desired symmetry type graph. In this paper we shall
show that some of the maniplexes constructed in~\cite{2OrbMani} are in
fact polytopal. More specifically, we shall see that if
$\{0,1,\ldots,n-1\}\setminus I$ has exactly two elements, there are
choices for the maniplex $\M$ and the monodromy $\eta$ that ensure
that the maniplex constructed in~\cite{2OrbMani} with symmetry type
$2^n_I$ is
polytopal. 

The first sections of this paper will introduce all the concepts that
we will use. In Section~\ref{s:basic} we define the basic concepts
such as \emph{abstract polytopes, flags, maniplexes, premaniplexes,
  automorphisms, symmetry type graphs, monodromies} and
\emph{coverings}. In Section~\ref{s:volts} we define the concepts of
\emph{fundamental groupoid} and \emph{fundamental group} of a graph,
\emph{voltage graphs, derived graphs}, and we enunciate
Theorem~\ref{t:IntProp}, which will be our main tool. In
Section~\ref{s:2alaM} we describe the construction $\hat{2}^\M$, which
gives rise to a family of regular polytopes of all ranks that we will
be using for all our examples. In Section~\ref{s:2OrbMani} we describe
the construction given in~\cite{2OrbMani} to get 2-orbit maniplexes,
and we also add some observations that will be helpful when proving
the main result. Finally, in Section~\ref{s:Poli2} we use this
construction and Theorem~\ref{t:IntProp} to prove that if $I$ has
exactly $n-2$ elements, then there exists a polytope $\p$ with
symmetry type $2^n_I$.

\section{Basic concepts}\label{s:basic}

An \emph{abstract polytope} is defined as a poset $(\p,\leq)$ that is
flagged, strongly flag-connected and satisfies the diamond
property. We clarify each of these properties now:

A poset is \emph{flagged} if it has a least element $F_{-1}$, a
greatest element $F_n$, and all maximal chains, called \emph{flags},
have the same finite size. If the size of the flags is $n+2$, we say
that the flagged poset has \emph{rank $n$}. Moreover, given an element
$F$ of a flagged poset, if the maximal chains having $F$ as their
biggest element have $r+2$ elements, we say that $F$ has \emph{rank
  $r$} and write $\rk(F)=r$. For the purposes of this paper, an
element of rank $i$ will be called an \emph{$i$-face}, and in general,
the elements of a flagged poset will be called \emph{faces}. The $(n-1)$-faces in a flagged poset of rank $n$ are called
\emph{facets}. If $\Phi$ is a flag, we denote its face of rank $i$ by
$(\Phi)_i$.

We say that a flagged poset satisfies the \emph{diamond condition} if
given two faces $F,G$ such that $F<G$ and $\rk(G)-\rk(F)=2$, there
exist exactly two faces $H_1,H_2$ such that $F<H_i<G$ for
$i\in\{1,2\}$.  We say that two flags $\Phi$ and $\Psi$ in a flagged
poset are \emph{$i$-adjacent} if they have the same faces with ranks
different than $i$ but they have different $i$-faces. If a flagged
poset satisfies the diamond condition, every flag $\Phi$ has exactly
one $i$-adjacent flag for $i\in\{0,1,\ldots,n-1\}$. The unique flag
$i$-adjacent to $\Phi$ is denoted $\Phi^i$.

We say that a flagged poset is \emph{strongly flag-connected} if given
two flags $\Phi, \Psi$, there is a sequence
$\Phi_0=\Phi, \Phi_1, \Phi_2, \ldots, \Phi_k=\Psi$, such that
$\Phi_{j-1}$ and $\Phi_j$ are $i_j$-adjacent for some $i_j$ satisfying
that $\Phi$ and $\Psi$ have different $i_j$-faces.

From now on, we will refer to abstract polytopes just as
\emph{polytopes}. We will also abuse notation and call $\p$ an
abstract polytope assuming that the symbol $\leq$ will be used for the
order relation (and $<$ for the strict order relation). A polytope of
rank $n$ is called an \emph{$n$-polytope}.

Given a polytope $\p$, the \emph{flag-graph of $\p$}, denoted
$\G(\p)$, is the edge-colored graph whose vertices are the flags of
$\p$ and $i$-adjacent flags are joined by an edge of color $i$. The
flag-graph of an $n$-polytope is an example of what is called an
\emph{$n$-maniplex} or a \emph{maniplex of rank $n$}. That is, a
connected simple graph with edges colored with the numbers
$\{0,1,\ldots,n-1\}$, in such a way that:
\begin{itemize}
\item every vertex is incident to an edge of each color, and
\item given two colors $i$ and $j$ such that $|i-j|>1$, the
  alternating paths of length 4 using edges of these two colors are
  closed.
\end{itemize}
The vertices of a maniplex are usually referred to as \emph{flags}.

At this point, we point out that although infinite maniplexes and
polytopes can be very interesting, for the purposes of this paper, all
maniplexes and polytopes are assumed to be finite. Some of the results
and definitions may also apply for the infinite case, but not all of
them (particularly, the main definition of Section~\ref{s:2alaM} needs
some fine-tuning for the infinite case).

Given a set of colors $I\subset\{0,1,\ldots,n-1\}$, and an
$n$-maniplex $\M$, we use the notation $\M_I$ for the subgraph of $\M$
induced by the edges with colors in $I$. We will denote the complement
of $I$ in $\{0,1,\ldots,n-1\}$ by $\ol{I}$, and if $I=\{i\}$, we denote
its complement just by $\ol{i}$. The connected components of
$\M_{\ol{i}}$ are the \emph{$i$-faces of $\M$}.

We say that an $i$-face $F$ and a $j$-face $G$ of a maniplex are
\emph{incident} if they have non-empty intersection. If in addition
$i\leq j$, we write $F\leq G$.

In this way we have defined a poset $\p(\M)$. If $\p$ is a polytope,
then $\p(\G(\p))$ is isomorphic to $\p$. On the other hand,
in~\cite{PolyMani} it is proven that if $\M$ satisfies certain
conditions called \emph{path intersection properties}, then $\p(\M)$
is a polytope and $\G(\p(\M))$ is isomorphic to $\M$, showing that
such maniplexes, called \emph{polytopal}, are exactly the flag-graphs
of polytopes.

An \emph{isomorphism} of $n$-maniplexes is just a graph isomorphism
that preserves the colors of the edges. Naturally, an automorphism of
a maniplex is just an isomorphism onto itself. The automorphism group
of a maniplex $\M$ is denoted by $\Gamma(\M)$. It is known that the
automorphism group of a maniplex acts freely on its flags (the proof
is identical to that of~\cite[Proposition~2A4]{ARP}, which is the
particular case of abstract polytopes). Moreover, the automorphisms of
the flag-graph of a polytope correspond to the automorphisms of the
polytope itself (as a poset) (see~\cite[Lemma~2A3]{ARP} for one
inclusion. The other inclusion can be proved using the construction
in~\cite[Section~3]{PolyMani}).

The \emph{symmetry type graph (STG)} of a maniplex $\M$, denoted by
$\T(\M)$, is simply the quotient of $\M$ by its automorphism
group. That is, the vertices of $\T(\M)$ are the flag orbits of $\M$
under its automorphism group, and two orbits are connected by an edge
of color $i$ if there is a pair of $i$-adjacent flags, one in each of
those orbits. If two $i$-adjacent flags are on the same orbit, we draw
a semi-edge of color $i$ on the vertex corresponding to that
orbit. Notice that $\T(\M)$ is not necessarily simple, but it
satisfies the other conditions that define an $n$-maniplex: it is
connected, each vertex is incident to exactly one edge of each color
in $\{0,1,\ldots,n-1\}$,
and the alternating paths of length 4 using two non-consecutive colors
are closed. We call such a graph an \emph{$n$-premaniplex}, or just a
\emph{premaniplex} if the rank $n$ is implicit.

The \emph{symmetry type graph} of a polytope $\p$ is just the
symmetry type graph of its flag-graph, and it is denoted by $\T(\p)$.

A \emph{$k$-orbit polytope (or maniplex)} is one with exactly $k$ flag
orbits under the action of its automorphism group. In other words, one
whose symmetry type graph has exactly $k$ vertices. A 1-orbit polytope
(or maniplex) is called \emph{regular}.

Given a maniplex $\M$, we define $r_i$ as the flag permutation that
maps each flag to its $i$-adjacent flag. It is easy to see that $r_i$
is an involution and that, if $|i-j|>1$, then $r_ir_j$ is also an
involution. In particular $r_i$ and $r_j$ commute whenever $i$ and $j$
are not consecutive.

The group $\Mon(\M)=\gen{r_i:0\leq i \leq n-1}$ is called the
\emph{monodromy group of $\M$}, and we will call each of its elements
a \emph{monodromy}. It is worth mentioning that some authors prefer
the term \emph{connection group} and the notation $\mathrm{Con}(\M)$. An
alternative way of defining the automorphism group of a maniplex $\M$
is as the set of flag permutations that commute with the elements of
the monodromy group.

For this paper, both the automorphism and the monodromy group will act
on the right.

If $\M$ is a regular maniplex (sometimes called \emph{reflexible}),
for a fixed flag $\Phi$ we can define $\rho_i$ as the (unique)
automorphism that maps $\Phi$ to $\Phi^i$. It is known that
$\{\rho_i:i\in\{0,1,\ldots,n-1\}\}$ is a generating set for the
automorphism group of $\M$ (see~\cite[Theorem 2B8]{ARP} for the
polytopal case). Moreover, the function mapping $\rho_i$ to $r_i$ can
be extended to a group anti-isomorphism between $\Gamma(\M)$ and
$\Mon(\M)$ (see, for example,~\cite[Theorem~3.9]{MixAndMon}
or~\cite[Section~7]{Maniplexes}). In particular, the monodromy group
of a regular maniplex acts regularly on its flags.

Given two $n$-premaniplexes $\M$ and $\X$, a \emph{covering
  projection} (\emph{covering} for short) from $\M$ to $\X$, is a
function $p$ from the vertices of $\M$ to the vertices of $\X$ that
preserves $i$-adjacencies. It can be proved that all coverings are
surjective. If $e$ is an edge $\M$ connecting vertices $x$ and $y$, we
define $p(e)$ as the edge connecting $p(x)$ and $p(y)$ that has the
same color as $e$. If there is a covering from $\M$ to $\X$, we say
that \emph{$\M$ covers $\X$}. Note that if $\M$ is a maniplex and $\X$
is its symmetry type graph, the natural projection from $\M$ to $\X$
is always a covering.

\section{Voltage graphs and intersection properties}\label{s:volts}

A \emph{dart} in a graph is just a directed edge. The \emph{inverse}
of a dart $d$, denoted by $d^{-1}$ is the dart that corresponds to the
same edge with the opposite orientation. We think of \emph{semi-edges}
as having only one orientation, having thus only one dart which is
inverse to itself. An edge with two different endpoints (and
therefore, with two different darts) will be called a
\emph{link}. Every dart $d$ has an \emph{initial vertex}, or
\emph{start-point} $I(d)$, and a \emph{terminal vertex} or
\emph{endpoint} $T(d)$, having the property that its underlying edge
is incident to both $I(d)$ and $T(d)$, and that $I(d^{-1}) = T(d)$ and
$T(d^{-1}) = I(d)$. For more on graphs with semi-edges
see~\cite{Voltajes}, for example.

A \emph{path} is a sequence of darts $W=d_1d_2\cdots d_k$, such that
$T(d_i) = I(d_{i+1})$ for every $i<k$. We say that $W$ \emph{goes from
  $u$ (its start-point) to $v$ (its endpoint)} if $u=I(d_1)$ and
$v=T(d_k)$. We say that a path is \emph{closed} if its start-point and endpoint are
equal. Otherwise, it is \emph{open}.

We also consider a formal empty path for each vertex that goes from
that vertex to itself. If $W$ is the empty path from $v$ to $v$, we
abuse notation and denote it by $v$ as well.

Two paths $W$ and $W'$ with the same start-point and endpoint are said
to be \emph{homotopic} if one can transform $W$ into $W'$ by a finite
sequence of the following operations:
\begin{itemize}
    \item Inserting two consecutive inverse darts at any point, that is $$d_1d_2\cdots d_i d_{i+1}\cdots d_k\mapsto d_1\cdots d_i d d^{-1} d_{i+1}\cdots d_k,$$ where $I(d)=T(d_i)$;
    \item Deleting two consecutive inverse darts at any point, that is $$d_1\cdots d_i d d^{-1} d_{i+1}\cdots d_k \mapsto d_1d_2\cdots d_i d_{i+1}\cdots d_k;$$
\end{itemize}
In this case we write $W\sim W'$.

It is known that homotopy is an equivalence relation and that it is
preserved by the usual concatenation of paths. We often abuse notation
and use $W$ as the name for a path or for its homotopy
class. Furthermore, even if we refer to $W$ as a path, we will always
be thinking of it up to homotopy.

The set of all homotopy classes of paths in a graph $X$, together with
the concatenation operation forms a groupoid called \emph{the
  fundamental groupoid of $X$}, which we will denote by $\Pi(X)$. The
subset of closed paths based at a vertex $x$ forms a group denoted by
$\Pi^x(X)$ and called \emph{the fundamental group of $X$ based at
  $x$}. It is easy to show that if $X$ is connected, all its
fundamental groups are isomorphic (moreover, they are conjugates by
elements of the fundamental groupoid).

Given a spanning tree $T$ and a vertex of a graph $X$ we can find a
distinguished set of generators for the fundamental group
$\Pi^x(X)$. For each dart $d$ in $X$ but not in $T$, take the path
$C_d$ that goes from $x$ to the starting point of $d$ through $T$,
then takes $d$, and then goes back to $x$ from the endpoint of $d$
trough $T$. It is easy to see that $\Pi^x(X)$ is generated by
$\{C_d\}$, where $d$ runs among the darts in $X$ not in $T$.

Given a graph $X$ and a group $\Gamma$, a \emph{voltage assignment (with
  voltage group $\Gamma$)} is a groupoid anti-morphism $\xi:\Pi(X)\to
\Gamma$. The pair $(X,\xi)$ is called a \emph{voltage graph}.

To construct a voltage assignment, it is enough to define it for the
darts of $X$ (in such a way that inverse darts have inverse
voltages). Then, the voltage of the path $W=d_1d_2,\cdots d_k$ is
simply $\xi(d_k)\cdots \xi(d_2)\xi(d_1)$.

Note that $\xi(\Pi^x(X))=\gen{\{\xi(C_d)\}_d}$, where $d$ runs among the
darts of the graph $X$ not in a given spanning tree $T$. This will be
important later, since very often we want to know the voltages of some
sets of closed paths. Knowing the voltages of a fundamental
group will also help us calculate the set of voltages of open paths
with fixed start-point and endpoint as a coset of the voltages of a
fundamental group.

Given a voltage graph $(X,\xi)$ with voltage group $\Gamma$, we can
construct the \emph{derived graph} $X^\xi$ as follows:
\begin{itemize}
\item The vertex set is $V\times \Gamma$ where $V$ is the vertex set of
  $X$.
\item The dart set is $D\times \Gamma$ where $D$ is the dart set of
  $X$.
\item If the dart $d$ goes from $x$ to $y$, the dart $(d,\gamma)$ goes
  from the vertex $(x,\gamma)$ to $(y,\xi(d)\gamma)$.
\item (Optional) If a dart (or vertex) $d$ has a color $c$, then the
  dart (or vertex) $(d,\gamma)$ also has color $c$.
\end{itemize}

Note that the inverse of the dart $(d,\gamma)$ is the dart
$(d^{-1},\xi(d)\gamma)$.

The voltage group $\Gamma$ always acts by automorphisms on the derived
graph $X^\xi$, with the action given by
$(x,\gamma)\sigma = (x,\gamma\sigma)$, where $x$ is a dart or vertex
of $X$ and $\gamma,\sigma\in \Gamma$. It should also be clear that if
$X$ has some coloring of its vertices or darts, then $\Gamma$
preserves the induced coloring when acting on $X^\xi$.

In~\cite{IntPropYo}, Hubard and the author of this paper prove the
following results:

\begin{lema}\label{l:DerivadaEsMani}
  {\rm\cite[Lemma~3.1]{IntPropYo}} Let $\X$ be a premaniplex and
  $\xi:\Pi(\X)\to \Gamma$ a voltage assignment such that $X$ has a
  spanning tree with trivial voltage on all its darts. Then $\X^\xi$
  is a maniplex if and only if the following conditions hold:
  \begin{itemize}
  \item The set $\xi(D)$ generates $\Gamma$ where $D$ is the set of
    darts of $\X$.
  \item If $d$ is the dart of a semi-edge, then $\xi(d)$ has order
    exactly two.
  \item If $d$ and $d'$ share both their start-point and their
    endpoint, then $\xi(d)\neq \xi(d')$.
  \item If $|i-j|>1$ every (closed) path $W$ of length 4 whose darts
    alternate between these two colors has trivial voltage.
  \end{itemize}
\end{lema}

Before giving the next result we need to introduce some notation.
Given two integers $k$ and $m$ we define $[k,m]:=\{k,k+1,\ldots,
m\}$. If $k>m$, then we define $[k,m]$ as the empty set.  If $\X$ is an
$n$-premaniplex, $a,b$ are two vertices in $\X$, and
$I\subset [0,n-1]$, then we define $\Pi^{a,b}_I(\X)$ as the subset of
$\Pi(\X)$ consisting of all the paths from $a$ to $b$ that only use
darts with colors in $I$. In particular, $\Pi^{a,b}_{[k,m]}(\X)$ is the
set of paths from $a$ to $b$ that only use colors between $k$ and $m$.

\begin{teo}\label{t:IntProp}
  {\rm\cite[Theorem~4.2]{IntPropYo}} Let $\X$ be a premaniplex of rank
  $n$ and $\xi:\Pi(\X)\to \Gamma$ a voltage assignment. Then $\X^\xi$
  is polytopal if and only if it is a maniplex and the following
  equation holds for all $k,m\in\{0,1,\ldots,n-1\}$ and all pairs of
  vertices $a,b\in \X$:
  \begin{equation}
    \label{eq:IntProp}
    \xi(\Pi^{a,b}_{[0,m]}(\X)) \cap \xi(\Pi^{a,b}_{[k,n-1]}(\X)) =
    \xi(\Pi^{a,b}_{[k,m]}(\X)).
  \end{equation}
\end{teo}

In this paper we will use Theorem~\ref{t:IntProp} to show that we can
construct 2-orbit polytopes with some prescribed symmetry type graphs.

\section{The construction $\hat{2}^\M$}\label{s:2alaM}

The main theorem in~\cite{2OrbMani} (Theorem~3) states that every
premaniplex with exactly two vertices is the symmetry type graph of a
maniplex $\tilde{\M}$. To prove this, given a regular maniplex $\M$ of
rank $n$ with some conditions, and a premaniplex $\X$ of rank $n+1$
with 2 vertices, the authors define a voltage assignment $\xi$ on $\X$
that gives a 2-orbit maniplex $\os{\M}$ with symmetry type graph
$\X$. The voltage group is a group acting on the set
$\M_W\times \Z_{2\ell}$, where $\M_W$ consists of a specific half of
the flags of $\M$ and $\ell$ (called $k$ in the original paper) is a
natural number (usually very large).

To be more precise, the conditions on $\M$ are (equivalent to) the
following:

\begin{enumerate}\label{Conds}
\item $\M$ has to cover $\X_{\ol{n}}$ (the premaniplex of rank $n$
  obtained by deleting the edges of color $n$ from
  $\X$)\footnote{In~\cite{2OrbMani} this is called having a coloring
    \emph{consistent with $I$} where $X=2^{n+1}_I$.}.
\item There is an involutory monodromy $\eta$ in $\M$ that maps all
  the flags of any given facet to different facets.
\end{enumerate}

A concrete family of maniplexes satisfying these conditions is given
and its elements are called $\M_n$ where $n$ denotes the rank. This
family is constructed recursively and it does not depend on the choice
of $\X$. More concretely, $\M_2$ is the flag-graph of the square, and
$\M_{n+1}$ is constructed from $\M_n$ as $\hat{2}^{\M_n}$ (which we
will define shortly). In~\cite{2ala} it
is proved that if $\p$ is a regular polytope then $\hat{2}^{\G(\p)}$
is the flag-graph of a regular polytope, and we will see shortly a
proof of the fact that if $\M$ is regular $\hat{2}^\M$ is regular
too. This implies that the family $\{\M_n\}_{n\geq 2}$ consists of
flag-graphs of regular polytopes. Since there are known examples of
2-orbit polyhedra (rank 3) with any given symmetry type, we are only
concerned about the family $\{\M_n\}_{n\geq 3}$, which are the
maniplexes used to construct 2-orbit maniplexes of ranks 4 and higher.

The construction $\hat{2}^\M$ works as follows: Given a maniplex $\M$,
the flags of $\hat{2}^\M$ are $\F(\M)\times \Z_2^{\Fac(\M)}$, where
$\F(\M)$ is the set of flags of $\M$ and $\Fac(\M)$ is the set of
facets. We will think of $\Z_2^{\Fac(\M)}$ as the set of functions
from $\Fac(\M)$ to the cyclic group $\Z_2$. Then the adjacencies in
$\hat{2}^\M$ are defined by:
\begin{eqnarray*}
  (\Phi,x)^i&:=&(\Phi^i,x) \quad {\rm if}\ i<n,\\
  (\Phi,x)^n&:=&(\Phi,x+\chi_{\Fac(\Phi)}),
\end{eqnarray*}
where $\Fac(\Phi)$ denotes the facet of $\Phi$ and given a facet $F\in
\Fac(\M)$ the vector $\chi_F$ is the one associated with the
characteristic function of $F$, that is, the vector with 1 in the
coordinate corresponding to $F$ and 0 in every other one.

If we remove the edges of color $n$ from $\hat{2}^\M$ we get one
connected component for each vector $x\in \Z_2^{\Fac(\M)}$. Each
component consists of all flags of type $(\Phi,x)$ where $\Phi$ is a
flag of $\M$. In particular, every facet of $\hat{2}^\M$ is isomorphic
to $\M$.

Given a polytope $\p$ we can construct a polytope $\hat{2}^\p$ such
that the flag-graph of $\hat{2}^\p$ is the maniplex
$\hat{2}^{\G(\p)}$. This construction is the following:

For $-1\leq i \leq n$, denote by $F_i$ the set of $i$-faces of
$\p$. If $f\in \Z_2^{\Fac(\M)}$ the {\em support of $f$},
denoted by $\supp(f)$, is defined as the set of facets $F\in \Fac(\M)$
such that $f(F) = 1$. If we denote by $\hat{F}_i$ the set of $i$-faces
of $\hat{2}^\p$, then for $-1\leq i \leq n$ let
$\hat{F}_i:=F_i\times \Z_2^{\Fac(\M)}/\sim$ where $(F,x)\sim (F',x')$
if and only if $F=F'$ and for every facet $G\in \supp(x+x')$
we have that $F\leq G$. We then add a formal greatest face
$F_{n+1}$. Finally, the incidence relation on $\hat{2}^\p$ is given by
$(A,x)<(B,y)$ if and only if $A<B$ and $(A,x)\sim (A,y)$.

We should remark that the construction $\hat{2}^\p$ is actually
$(2^{\p^*})^*$, where the superscript $*$ denotes duality and $2^\p$
is a better-known construction (see~\cite{2alaDanzer}).

In~\cite{twist} it is mentioned that the operation $\hat{2}^\M$ where
$\M$ is a maniplex generalizes $\hat{2}^\p$ where $\p$ is a
polytope. What the authors mean by this is that if $\p$ is a polytope,
then $\G(\hat{2}^\p)$ is isomorphic to $\hat{2}^{\G(\p)}$
(see~\cite[Section~5.1]{MiTesis} for an elementary proof).

The automorphisms of $\M$ have a natural action on $\hat{2}^\M$. First
we define a left action of $\Aut(\M)$ on $\Z_2^{\Fac(\M)}$ as follows:
given an automorphism $\sigma$ of $\M$, a vector $x$ in
$\Z_2^{\Fac(\M)}$ and a facet $F$ of $\M$ we define that
$\sigma x:f\mapsto x(f\sigma)$. Now, if $\sigma$ is an automorphism of
$\M$, we can extended it to an automorphism of $\hat{2}^\M$ by
defining $(\Phi,x)\sigma := (\Phi\sigma, \sigma^{-1} x)$. Note that
$\supp(\sigma^{-1}x) = \supp(x)\sigma$, so we are actually using the
natural action of an automorphism on a set of facets.

In~\cite[Section 6]{twist} it is proved that $\Gamma(\M)$ this action
on $\hat{2}^\M$  is by automorphisms. Moreover, the
automorphism group of $\hat{2}^\M$ is a product of
$\Gamma(\M)$and the group
$T=\{T_y:(\Phi,x)\mapsto (\Phi,x+y) : y\in \Z_2^{\Fac(\M)}\}$ . This
implies that if an $n$-maniplex $\M$ has symmetry type graph $\X$ then
the symmetry type graph of $\hat{2}^\M$ is obtained by adding
semi-edges of color $n$ to each vertex of $\X$. In particular, if $\M$
is regular, then $\hat{2}^\M$ is regular as well.

When we let $\M_2$ be the flag-graph of a square and define
$\M_{i+1}:=\hat{2}^{\M_i}$ for $i\geq 2$, $\M_3$ happens to be the
flag graph of the map on the torus called $\{4,4\}_{(4,0)}$ which can
be thought of as a $4\times 4$ chess board in which we identify
opposite sides (without twisting). The 3-faces of every subsequent
$\M_n$ will be of this type and we will make use of this in our
proofs.

A {\em lattice} is a poset in which every pair of elements
$\{A,B\}$ has a lowest upper bound $A\vee B$ (called the {\em join}
of $A$ and $B$) and a greatest lower bound $A\wedge B$ (called the
{\em meet} of $A$ and $B$). It is known that if $\p$ is a lattice then
$\hat{2}^\p$ is also a lattice. One can verify that if $A$ and $B$ are
faces of $\p$, $x,y \in \Z_2^{\Fac(\p)}$ and there is at least one
facet of $\hat{2}^\p$ containing both $(A,x)$ and $(B,y)$ then,
\[
  (A,x) \vee (B,y) = (A\vee B,x)\sim (A\vee B, y),
\]
and
\[
  (A,x)\wedge (B,y) = (C,x) \sim (C,y)
\]
where $C$ is the greatest lower bound of the set $\supp(x+y) \cup
\{A,B\}$ . For a full proof of a more general result
see~\cite[Theorem~5]{Egon2alaLattice}.

\begin{coro}\label{c:lattice}
  For $n\geq 2$, the poset associated with $\M_n$ is a lattice.
\end{coro}

In~\cite[Proposition~15,~Lemma~17]{2OrbMani} it is proved that the
family $\{\M_n\}_{n\geq 3}$ satisfies both conditions mentioned at the beginning of this section. To prove that it satisfies
the second condition, the authors prove and use the following lemma:

\begin{lema}\label{l:EtaExiste}{\rm\cite[Lemma~16]{2OrbMani}}
  If $\M$ has a set of facets $S$ which is not invariant under any
  non-trivial automorphism, then there is an involutory monodromy
  $\eta$ of $\hat{2}^\M$ that maps all the flags of any given facet to
  different facets.
\end{lema}

Such a set $S$ is constructed for $\M_3$, and recursively constructed
for $\M_{n+1}$ in terms of one constructed for $\M_n$. The
construction is as follows:

If $S\subset \Fac(\M)$ is a set of facets that is not invariant under
any non-trivial automorphism, then
$\hat{S}:=\{(F_n,\chi_F)|F\in S\}\cup\{(F_n,0)\}$ (here 0 denotes the
all-zero vector in $\Z_2^{\Fac{\M}}$) is a set of facets of
$\hat{2}^\M$ and it satisfies the same condition.

In~\cite{2OrbMani} the authors only care that $S$ is not invariant
under any non-trivial automorphism, however we shall choose an $S$
that satisfies some extra conditions, which will prove useful when
dealing with the polytopality of the constructed 2-orbit
maniplexes.

\begin{figure}
  \begin{center}
    \includegraphics[width=5cm]{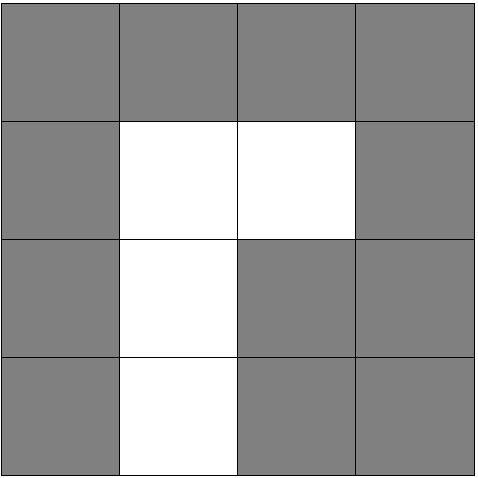}
    \caption{\label{f:S3} The set $S_3$ consisting of the shaded
      2-faces is not invariant under any non-trivial symmetry of
      $\M_3$ and it is not contained in the closure of two 0-faces.}
  \end{center}
\end{figure}

From now on, we call the shaded set in Figure~\ref{f:S3} $S_3$, and we
call $S_n\subset\Fac(\M_n)$ the set constructed recursively as
$S_n:=\hat{S}_{n-1}$ for $n\geq 4$. Note that $S_3$ is not invariant
under non-trivial automorphisms of $\M_3$, and therefore $S_n$ is not
invariant under non-trivial automorphisms of
$\M_n$. In~\cite{2OrbMani}, the authors actually use the complement of
the set we have chosen as $S_3$. The reason why we have chosen the
complement will be apparent when we prove Lemma~\ref{l:DosCaras}.

\begin{Def}\label{d:Cerradura}
  Given a polytope $\p$ and a face $f\in \p$ let us define the {\em
    closure of $f$}, denoted by $\ol{f}$, as the set of all the facets
  of $\p$ which are incident to $f$.
\end{Def}

Given a pre-ordered \footnote{A pre-order $\leq$ on a set $\p$ is a
                   relation that is transitive and reflexive, but it
                   may not be anti-symmetric.}
set $(\p,\leq)$, one can give it a topology by defining that a subset
$C$ is closed if and only if whenever $x\in C$ and $x\leq y$ we get
that $y\in C$. In fact every topology on a finite set can be obtained
from a pre-order in this way (see, for example,~\cite{ETF2}). In such
context, the closure of a set $S$ is just
$\ol{S}=\{y\in \p: \exists x\in S, x\leq y\}$. If we remove the least
and greatest faces of a polytope $\p$ and equip it with this topology,
then what we are calling the closure of $f$ is actually just the set
of facets contained in the topological closure of $\{f\}$.

\begin{lema}\label{l:DosCaras}
  Given any two proper faces $u,v$ of $\M_n$ with $n\geq 3$, the set $S_n$ is not
  contained in $\ol{u} \cup \ol{v}$.
\end{lema}
\begin{dem}
  We proceed by induction on $n$. For $n=3$ this is a simple
  observation obtained from Figure~\ref{f:S3}. Note that the closure
  of a face is contained in the closure of any incident face of
  smaller rank, so we only need to prove the lemma for
  0-faces. Suppose the lemma is true for $\M_n$. Let $u,v$ be 0-faces
  of $\M_{n+1}$, so that $u=(u',x)$ and $v=(v',y)$ for some 0-faces
  $u',v'$ of $\M_n$ and $x,y \in \Z_2^{\Fac(\M_n)}$. We will prove
  that $\ol{u} \cup \ol{v}$ does not contain even
  $S_{n+1}\setminus \{(F_n,0)\}$, so we may assume that both $u$ and
  $v$ are each incident to at least one element of $S_{n+1}$ other
  than $(F_n,0)$, say for example that $u=(u',x)<(F_n,\chi_G)$ for
  some facet $G$ of $\M_n$. By our definition of the order $<$, this
  implies that $u'<F_n$ (which is tautological) and
  $(u',x)\sim(u',\chi_G)$, so we may assume without loss of generality
  that $x = \chi_G$. Analogously we may assume that there exists a
  facet $F$ in $\M_n$ such that $y =\chi_F$.

  A facet $(F_n,z)$ is in $\ol{u} \cup \ol{v}$ if and only if
  $\supp(z+\chi_G)\subset \ol{u'}$ or
  $\supp(z+\chi_F)\subset \ol{v'}$. By induction hypothesis, $S_n$ is
  not contained in $\ol{u'} \cup \ol{v'}$, so there exists a facet
  $D \in S_n \setminus (\ol{u'} \cup \ol{v'})$. Then $(F_n,\chi_D)$
  is a facet in $S_{n+1}$ not in $\ol{u} \cup \ol{v}$.
  \qed \end{dem}

In topological terms, Lemma~\ref{l:DosCaras} tells us that no set of
two elements is {\em dense} in $S_n$.

Now we turn our attention to the monodromy $\eta$. In the proof of
Lemma~\ref{l:EtaExiste} found in~\cite{2OrbMani}, $\eta$ is
constructed as follows:

Let $\M$ be a regular maniplex and let $S$ be a set of facets not
invariant under non-trivial automorphisms. For every facet $F \in S$
let $\Phi_F$ be a fixed base flag in that facet. Let $F_0$ be a base
facet in $S$ and let $\Phi = \Phi_{F_0}$. For every $F\in S$, let
$\omega_F$ be a monodromy of $\M$ that maps $\Phi$ to $\Phi_F$. Note
that since $\M$ is regular, its monodromy group acts
regularly on its flags, so $\omega_F$ is actually unique. Then for
every flag $(\Psi,x)$ in $\hat{2}^\M$ define
$(\Psi,x)\eta:=(\Psi,x+\sum_{F\in S}\chi_{\Fac(\Psi\omega_F)})$. The
action of $\eta$ will be more clear with the following lemma.

\begin{lema}\label{l:AccionEta}
  Let $(\Psi,x) \in \hat{2}^\M$. Let $\gamma$ be the automorphism of
  $\M$ mapping the base flag $\Phi$ of $F_0$ to $\Psi$. Then, the
  vector corresponding to the facet of $(\Psi,x)\eta$ differs from $x$
  only in the coordinates corresponding to $S\gamma$, that is,
  if $(\Psi,x)\eta = (\Psi,y)$ then $\supp(x+y) = S\gamma$.
\end{lema}

\begin{dem}
  For every $F \in S$ we have
  that
  \[
  \Psi\omega_F = (\Phi\gamma)\omega_F = (\Phi\omega_F)\gamma,
  \]
  so
  \[
  \Fac(\Psi\omega_F) = \Fac((\Phi\omega_F)\gamma) =
  (\Fac(\Phi\omega_F))\gamma =  (\Fac(\Phi_F))\gamma = F\gamma.
  \]

  This implies that
  $(\Psi,x)\eta = (\Psi,x+\sum_{F\in S} \chi_{\Fac(\Psi\omega_F)}) =
  (\Psi,x+\sum_{F\in S} \chi_{F\gamma})$ so
  $y = x+\sum_{F\in S}\chi_{F\gamma}$. By doing the change of variable
  $G=F\gamma$ we get $y=x+\sum_{G\in S\gamma}\chi_G$ and
  $\supp(x+y) = S\gamma$. \qed
\end{dem}

Given two sets $S$ and $S'$ of facets of a polytope $\p$, we say that
$S'$ is {\em a copy of $S$} if there is an automorphism $\gamma$ of
$\p$ such that $S'=S\gamma$. Lemma~\ref{l:AccionEta} tells us that if
$(\Psi,x)\eta = (\Psi,y)$ then $\supp(x+y)$ is a copy of $S$.

In~\cite[Lemma~16]{2OrbMani} it is proved that $\eta$ is in fact a
monodromy of $\hat{2}^\M$ and that if $S$ is not invariant under
non-trivial automorphisms of $\M$ then $\eta$ maps all the flags of a
facet of $\hat{2}^\M$ to different facets.

Note that using a set of facets of $\M_n$ not invariant under
non-trivial automorphisms we constructed the monodromy $\eta$ for
$\M_{n+1}$. We have found such sets of facets for $\M_n$ with
$n\geq 3$, which means that we have found the monodromy $\eta$ for
$n\geq 4$. We have not found the monodromy $\eta$ for $\M_3$, since
every set of facets (edges) of $\M_2$ (the square) is invariant under
some non-trivial automorphism. Actually $\M_3$ does not have such a
monodromy. 

However, the map on the torus $\{4,4\}_{(8,0)}$ (a chess board of
regular size where each border is identified with its opposite) does
admit such a monodromy: simply take $\eta=r_2r_1r_0r_1r_2r_1r_2r_1$.
This monodromy is indeed involutory (use the fact that $r_1r_2$ has
order 4) and it maps the 8 flags of a facet to the 8 squares where a
knight could legally move in a game of chess. This example also covers
all 2-vertex premaniplexes in rank 3, so, to construct 2-orbit
polytopes of rank 4, one may use this example instead of $\M_3$.

\section{2-orbit maniplexes}\label{s:2OrbMani}

The symmetry type of a 2-orbit maniplex is denoted by $2_I^n$ where
$n$ denotes the rank and $I$ is a proper subset of $\{0,1,\ldots, n-1\}$.  Some authors
assume the rank is implicit and just write $2_I$, but we shall not do
that in this paper as we will be working with several different ranks at
once. In this notation, $I$ is the set of colors $i$ such that any
flag is in the same orbit as its $i$-adjacent flag. In terms of the
symmetry type graph, $I$ is the set of colors of the semi-edges.
Informally, we may think of $I$ as the set of ``allowed kinds of
reflections''. In particular, chiral $n$-polytopes are those of type
$2_\emptyset^n$, that is, those in which no reflections are allowed.

In order to find 2-orbit polytopes (resp. maniplexes) with every
possible 2-vertex premaniplex as its symmetry type graph, we only need
to find examples with all the symmetry types $2^n_I$ where
$0,n-1 \in \ol{I}$, and then use the operations $\hat{2}^\M$ and
duality (see~\cite[Corollary~12]{2OrbMani} for details). So without loss of
generality, we will assume that $0,n-1 \in \ol{I}$.

Given a 2-vertex premaniplex $\X = 2_I^{n+1}$ with $0,n\in \ol{I}$,
the authors of \cite{2OrbMani} give a voltage assignment to its darts
such that the derived graph is a 2-orbit maniplex with symmetry type
graph isomorphic to $\X$. We discuss this voltage assignment next.

Let us color the vertices of $\X$ one white and one black. Let
$\{r_0,r_1,\ldots,r_{n-1}\}$ be the distinguished generators of the
monodromy group of $\M_n$. If $d$ is a semi-edge incident on the white vertex we assign to it the
voltage $r_i$ where $i$ is its color. If $d$ is a semi-edge of color
$i$ incident to the black vertex, we assign to it the voltage
$r_0r_ir_0$, or in other words, $r_i$ if $i>1$ and $r_0r_1r_0$ if
$i=1$. Finally, if $d$ is a dart of color $i<n$ from the white vertex
to the black vertex we assign to it the voltage $r_0r_i$, in
particular, the voltage of the edge of color 0 is trivial. We shall
give the voltage of the edge of color $n$ shortly (which will be an
involution, so orientation is irrelevant). From now on, we call this
voltage assignment $\xi$.

\begin{figure}
  \begin{center}
    \includegraphics[width=5cm]{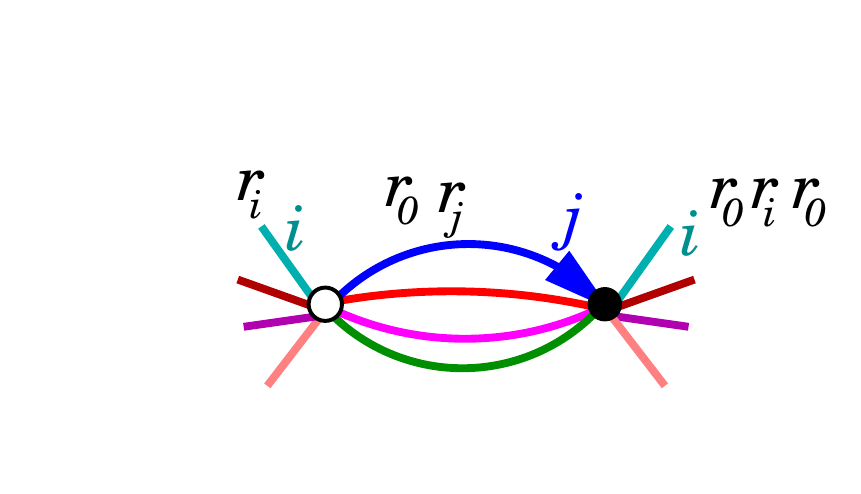}
    \caption{\label{f:Volts2orb} The voltage $\xi$ we use on
      $\X=2^{n+1}_I$.}
  \end{center}
\end{figure}

Let $p:\M_n\to \X_{\ol{n}}$ be a covering (remember that $\X_{\ol{n}}$
is the result of removing the edge of color $n$ from $\X$). We color
the flags of $\M_n$ white or black according to the color of their
image under $p$. Since $\X_{\ol{n}}$ has exactly 2 vertices and $\M_n$
covers $\X$, then every monodromy of $\M_n$ either preserves the color
of every flag or it changes the color of every flag. Note that the
voltages of every dart preserve the color of the flags when thought of
as monodromies.

\subsection{The voltage of the darts of color $n$}\label{s:y_n}

Now we turn our attention to the voltage of the darts of color $n$. In \cite{2OrbMani} the authors make the voltage group act on
$(\M_n)_w \times \Z_{2\ell}$ for some large $\ell$, where $(\M_n)_w$ is the
set of white flags of $\M_n$ (recall that all voltages preserve the
color of the flags). This action is simply defined by
$(\Psi,x)\omega = (\Psi\omega, x)$ for all
$(\Psi,x) \in (\M_n)_w \times \Z_{2\ell}$ and every monodromy
$\omega$. In other words, each monodromy acts as usual in the
$(\M_n)_w$-coordinate and as the identity on the
$\Z_{2\ell}$-coordinate. The voltage $\os{y}_n$ of the darts of color $n$
is defined as the composition of three commuting involutions
$\os{\rho}_0$, $r_0$ and $s$. Both $r_0$ and $\os{\rho}_0$ act only on
the $(\M_n)_w$ coordinate and are independent of the $\Z_{2\ell}$ one,
while $s$ acts only on the $\Z_{2\ell}$ coordinate. To avoid future
confusion, the reader must recall that
$\F(\M_n) = \F(\M_{n-1}) \times \Z_2^{\Fac(\M_{n-1})}$, so
$\F(\M_n)\times \Z_{2\ell}$ may be regarded as
$\F(\M_{n-1}) \times \Z_2^{\Fac(\M_{n-1})} \times \Z_{2\ell}$. The
importance of the $\Z_{2\ell}$ coordinate lies in that it is used to
prove that the derived maniplex is in fact a 2-orbit maniplex (and not
a regular one), but it has no relevance when proving that the derived
graph is a maniplex. For our purposes, we will ignore $s$ for
now. Therefore, we may think of all actions as being on
$\F(\M_n) = \F(\M_{n-1}) \times \Z_2^{\Fac(\M_{n-1})}$, and assign the
voltage $y_n := \os{\rho}_0r_0$ to the edge of color $n$.

In order to define $\os{\rho}_0$ we first need to choose a {\em base
  flag} $\Phi_F$ for each facet $F$ of $\M_n$. For the proof in
\cite{2OrbMani} to work, we first choose a {\em base facet}
$F_0=(f_{n-1},0)$ of $\M_n$ (where $f_{n-1}$ is the greatest face of
$\M_{n-1}$) and then for every {\em white} flag $\Psi$ in $F_0$ we set
$\Psi r_0 \eta r_0$ to be the base flag of its facet. Recall that
$\eta$ maps all the flags of $F_0$ to different facets, so there is no
ambiguity. In~\cite{2OrbMani} all the other base flags are chosen
arbitrarily, however we will later have a preferred choice for them
too. Then we define $\os{\rho}_0$ as the flag-permutation that acts on
each facet $F$ as the reflection (facet automorphism) that fixes all
the faces of the base flag $\Phi_F$ but its 0-face.

Note that if we replace a base flag with any flag sharing the same
edge (1-face) and the same facet, we get the same permutation
$\os{\rho}_0$. So we may define a {\em base edge} of a facet $F$ as
the 1-face of the corresponding base flag and forget about the base
flag.

Now the previous choice of base edges would be equivalent to the
following: For every {\em black} flag $\Psi$ in $F_0$ we set
$(\Psi\eta)_1$ (the edge of $\Psi\eta$) to be the base edge of its
facet.

Let $\Psi$ be a black flag in $F_0$. Since $\M_n = \hat{2}^{\M_{n-1}}$,
and by the choice of $F_0$, the flag $\Psi$ can be written as
$(\psi,0)$ for some flag $\psi$ in $\M_{n-1}$. Lemma~\ref{l:AccionEta} tells
us that $\Psi\eta = (\psi,x)$ for some $x\in\Z_2^{\Fac(\M_{n-1})}$
satisfying that $\supp(x)$ is a copy of
$S_{n-1}$.

So for facets corresponding to vectors whose support is a copy of
$S_{n-1}$ we are forced to choose a specific base edge, but for any
other facet we may choose the base edge as we want. Let $(e,0)$ be the
base edge of the base facet $F_0 = (f_{n-1},0)$. Then for every $x$
whose support is not a copy of $S_{n-1}$ we choose $(e,x)$ as the base
edge of the facet $(f_{n-1},x)$. We summarize this in the following
definition:

\begin{Def}\label{d:AristaBase}
  {\em Base edges:} Let $(e,0)$ be the base edge of the base facet of
  $\M_n$, let $x\in\Z_2^{\Fac(\M_n)}$ and let $f_{n-1}$ be the
  greatest face of $\M_{n-1}$. If $\supp(x)$ is not a copy of
  $S_{n-1}$, then we define the base edge of the facet $(f_{n-1},x)$
  of $\M_n$ to be $(e,x)$. If $\supp(x)=S_{n-1}\gamma$ for some
  $\gamma\in \Gamma(\M_{n-1})$, then we define the base edge of
  $(f_{n-1},x)$ to be $(e\gamma,x)$.
\end{Def}

\begin{coro}\label{c:AristaBase}
  If $x\in \Z_2^{\Fac(\M_n)}$ is such that $\supp(x)\subset \ol{u}
  \cup \ol{v}$ for some faces $u,v$ in $\M_{n-1}$, then the base edge
  of the facet $(f_{n-1},x)$ is $(e,x)$.
\end{coro}

\begin{coro}\label{c:AristaBaseBase}
  If $(e,0)$ is the base edge of the base facet $(f_{n-1},0)$ of
  $\M_n$, then it is also the base edge of any other facet containing
  it.
\end{coro}
\begin{dem}
  If $(f_{n-1},x)$ is a facet containing $(e,0)$ then $\supp(x)
  \subset \ol{e}$. Then, Corollary~\ref{c:AristaBase} tells us that
   base edge of $(f_{n-1},x)$ is $(e,x) \sim (e,0)$.\qed
\end{dem}

Finally we let $y_n:=\os{\rho}_0r_0$ and extend the voltage
assignment $\xi$ by assigning $y_n$ as the voltage of the edge of
color $n$ of $\X$. Note that since $\os{\rho}_0$ acts as an
automorphism in each facet, it commutes with all the monodromies
that do not use the generator $r_n$. In particular $\os{\rho}_0$
commutes with $r_0$, implying that $y_n$ is an involution.

As previously discussed, in~\cite{2OrbMani} the authors consider the
voltage group as acting on $(\M_n)_w\times \Z_{2\ell}$ where $(\M_n)_w$
is the set of white flags of $\M_n$ and $\ell$ is some large integer. If
$(\Psi,a) \in (\M_n)_w\times \Z_{2\ell}$ and $\omega$ is a color
preserving monodromy of $\M_n$ then $(\Psi,a)\omega$ is simply defined
as $(\Psi\omega,a)$. Then they use $\os{y}_n:=\os{\rho}_0r_0s$ as the
voltage of the edge of color $n$, where $s$ is an involution acting
only on the $\Z_{2\ell}$ coordinate. We will call the voltage assignment
used in~\cite{2OrbMani} $\xi'$. In~\cite{2OrbMani} the authors prove
that $\X^{\xi'}$ is a maniplex with symmetry type graph isomorphic to
$\X$. The proof of the fact that $\X^{\xi'}$ is a maniplex also
applies to $\X^\xi$, but the proof of the fact $\X^{\xi'}$ is not
regular relies heavily on $\ell$ being
large. On the other hand $\X^\xi$ is either regular or has symmetry type
graph isomorphic to $\X$. However our proofs to show polytopality of
some derived maniplexes will be much clearer if we work with the
voltage assignment $\xi$ and consider that $\xi=\pi_{\M_n}\xi'$ where
$\pi_{\M_n}$ is the projection to the $\M_n$-coordinate. That is, the
$\xi$-voltage of a path is just the first coordinate of the
$\xi'$-voltage of the same path. More formally, we define:

\begin{equation}
  \label{eq:xi}
  \xi(d)=
  \begin{cases}
    1  & \text{If }  d \text{ has color } 0\\
    r_i & \text{If } d \text{ is a semi-edge based on } a \text{ with color } i
    \in
    [1,n-1]\\
    r_0r_ir_0  & \text{If }  d \text{ is a semi-edge based on } b \text{ with color
    } i \in
    [1,n-1]\\
    r_0r_i  & \text{If }  d \text{ goes from } a \text{ to } b \text{ and has color
    }
    i \in [1,n-1]\\
    r_ir_0  & \text{If }  d \text{ goes from } b \text{ to } a \text{ and has color
    }
    i \in [1,n-1]\\
    y_n=\os{\rho}_0r_0  & \text{If }  d \text{ has color } n
  \end{cases}
\end{equation}
and

\begin{equation}
  \label{eq:xi'}
  \xi'(d)=
  \begin{cases}
    1  & \text{If }  d \text{ has color } 0\\
    r_i  & \text{If }  d \text{ is a semi-edge based on } a \text{ with color } i
    \in
    [1,n-1]\\
    r_0r_ir_0  & \text{If }  d \text{ is a semi-edge based on } b \text{ with color
    } i \in
    [1,n-1]\\
    r_0r_i  & \text{If }  d \text{ goes from } a \text{ to } b \text{ and has color
    }
    i \in [1,n-1]\\
    r_ir_0  & \text{If }  d \text{ goes from } b \text{ to } a \text{ and has color
    }
    i \in [1,n-1]\\
    \os{y}_n=\os{\rho}_0r_0s  & \text{If }  d \text{ has color } n
  \end{cases}
\end{equation}
Note that two paths may have the same
$\xi$-voltage while having different $\xi'$-voltages, but not the
other way around.

\section{Polytopality}\label{s:Poli2}

\begin{figure}
  \centering
  \includegraphics[width=5cm]{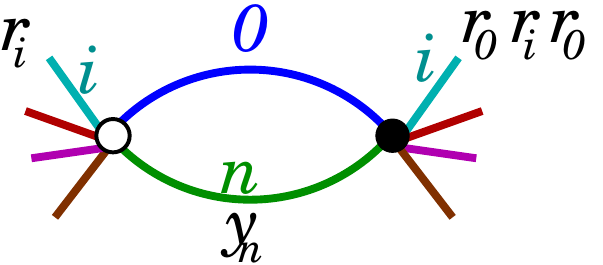}
  \caption{The premaniplex $2^{n+1}_{[1,n]}$ with its voltage
    assignment.}
  \label{f:2links}
\end{figure}

In this section we prove that if $I=[1,n-1]=\{1,2,\ldots,n-1\}$, then the
derived graph from the voltage graph $\X = 2^{n+1}_I$ in
Figure~\ref{f:Volts2orb} is polytopal.
We will first prove that most of
the intersection properties are satisfied for arbitrary $I$,
implying that there are only a few intersection properties that would
be needed to check to prove that there are 2-orbit polytopes of all
symmetry types in rank $n\geq 3$. Once we have done this, we will prove
that the remaining intersection properties are satisfied for
$I=[1,n-1]$ (see Figure~\ref{f:2links}).

In order to deal with the intersection properties in
Theorem~\ref{t:IntProp} for this particular voltage assignment, we
should first understand better what are the voltages of paths in the
voltage graph in Figure~\ref{f:Volts2orb} with respect to the colors
they use. To do this we prove Lemmas~\ref{l:VoltsCerr}
and~\ref{l:VoltsAbs}, which describe the voltages of some of the
paths relevant to our results. 

Let us first prove the following:

\begin{af}\label{a:Volts2}
  Let $n\geq 3$. Let
  $W=d_1d_2\cdots d_m$ be a path in the voltage graph of
  Figure~\ref{f:Volts2orb}. Let $c_i$ be the color of the dart
  $d_i$. If $c_i\notin \{1,n\}$ for all $i\in\{1,2,\ldots,m\}$ then
  $\xi(W)=r_0^\varepsilon r_{c_m}r_{c_{m-1}}\cdots r_{c_1}$ where
  $\varepsilon$ is 0 or 1 depending on whether the path $W$ is closed
  or open, respectively.
\end{af}
\begin{dem}
  Note that $W$ is closed if it uses an even number of
  links, and it is open if it uses an odd number of links.

  If $d_m$ is a semi-edge based on the white vertex, then
  \[\xi(d_1d_2 \cdots d_m) = r_{c_m}\xi(d_1d_1\cdots d_{m-1}).\] If
  $d_m$ is a dart from the white vertex to the black one, then
  \[\xi(d_1d_2 \cdots d_m) = r_0r_{c_m}\xi(d_1d_1\cdots d_{m-1}).\] If
  $d_m$ is a dart from the black vertex to the white one, then, since
  $c_m\neq n$,
  \[\xi(d_1d_2 \cdots d_m) = r_{c_m}r_0\xi(d_1d_1\cdots d_{m-1})\] but
  since $c_m\neq 1$ we know that $r_{c_m}r_0=r_0r_{c_m}$ and we get
  \[\xi(d_1d_2 \cdots d_m) = r_0r_{c_m}\xi(d_1d_1\cdots d_{m-1}).\] And
  if $d_m$ is a semi-edge on the black vertex then
  \[\xi(d_1d_2 \cdots d_m) = r_0r_{c_m}r_0\xi(d_1d_1\cdots d_{m-1})\]
  but since $c_m\neq 1$ this means
  $\xi(d_1d_2 \cdots d_m) = r_{c_m}\xi(d_1d_1\cdots d_{m-1})$. If we
  repeat this argument for $d_{m-1}$, then $d_{m-2}$ and so on, and
  note that since $c_i\neq 1$ then $r_{c_i}$ commutes with $r_0$ for
  all $i$, we get the desired result. \qed
\end{dem}

If $\Phi$ is a flag and $K\subset\{0,1,\ldots,n-1\}$, we will denote
by $(\Phi)_K$ the set of faces in $\Phi$ whose rank is in
$K$. The following lemma characterizes the voltages of {\em closed}
paths that do not use edges of color $1$ or the edge of color $n$ and
voltage $y_n$.

\begin{lema}\label{l:VoltsCerr}
  Let $\omega$ be a monodromy of $\M_n$ that preserves the color of
  its flags and $K\subset \{0,1,\ldots,n-1\}$. Suppose $1\in K$. Then
  for every white flag $\Phi$ of $\M_n$ we have that
  $(\Phi)_K = (\Phi\omega)_K$ if and only if $\omega$ is the voltage
  of a {\em closed} path based on the white vertex of $\X$ that does
  not use colors in $K\cup \{n\}$.
\end{lema}
\begin{dem}
  If $(\Phi)_K = (\Phi\omega)_K$, by strong connectedness of $\M_n$
  there is a path $\os{W}$ from $\Phi$ to $\Phi\omega$ not using
  colors in $K$. Let $c_1c_2\cdots c_k$ be the sequence of colors of
  $\os{W}$. Then $\Phi\omega = \Phi r_{c_1}r_{c_2}\cdots r_{c_k}$, but
  since $\M_n$ is regular, the action of the monodromy group on the
  flags is regular, so we conclude that
  $\omega = r_{c_1}r_{c_2}\cdots r_{c_k}$.

  Let $W$ be the path on $\X$ that starts on the white vertex and
  follows the sequence of colors $c_1c_2\cdots c_k$, that is
  $W = p(\os{W})$. We know that $\Phi\omega$ and $\Phi$ are both
  white, so $W$ must be a closed path. Recall that the voltage of $W$
  is the product of the voltages of its darts but in reverse
  order. Let $a_1a_2\cdots a_k$ be the dart sequence of $W$. Let us
  consider $W^{-1} = a_ka_{k-1}\cdots a_1$. Since $W$ does not use
  colors in $K$ and $1\in K$, Claim~\ref{a:Volts2} tells us that
  $\xi(W^{-1}) = r_0^\varepsilon\omega$ and since $W$ is closed,
  $W^{-1}$ is closed too, so $\varepsilon = 0$ and
  $\xi(W^{-1})=\omega$. Thus we have found a closed path that does not
  use colors in $K\cup\{n\}$ and has voltage $\omega$.

  For the converse, let $W=a_1a_2\cdots a_k$ be a closed path based on
  the white vertex of $\X$ and suppose that $W$ does not use colors in
  $K\cup \{n\}$. Since $W$ does not use the color 1 (because $1\in K$),
  Claim~\ref{a:Volts2} tells us that
  $\omega:=\xi(W)=r_{c_k}r_{c_{k-1}}\cdots r_{c_1}$ where $c_i$ is the
  color of $a_i$. Since $c_i\notin K\cup \{n\}$ we know that
  $(\Phi\omega)_K = (\Phi)_K$. \qed
\end{dem}

Using the same logic one can prove the following lemma, which
characterizes the voltages of {\em open} paths that do not use edges
of color $1$ or the edge of color $n$ and voltage $y_n$.

\begin{lema}\label{l:VoltsAbs}
  Let $\omega$ be a monodromy of $\M_n$ that preserves the color of
  the flags and $K\subset \{0,1,\ldots,n-1\}$. Then
  $(\Phi)_K = (\Phi^0\omega)_K$ if and only if $\omega$ is the voltage
  of an {\em open} path in $\X$ that does not use colors in
  $K \cup \{n\}$.
\end{lema}

We want to show that if $\X$ is the 2-vertex premaniplex with links of
colors 0 and $n$ and semi-edges of the colors in-between (see
Figure~\ref{f:2links}), and $\xi$ is the voltage assignment defined
in~(\ref{eq:xi}), then the condition of Theorem~\ref{t:IntProp} is
satisfied, that is, we want to prove that if
$\X = \hat{2}^{n+1}_{[1,n-1]}$, we have that

\[
\xi(\Pi^{a,b}_{[0,m]}(\X)) \cap \xi(\Pi^{a,b}_{[k,n]}(\X)) =
\xi(\Pi^{a,b}_{[k,m]}(\X)),
\]
for all $k,m \in [0,n]$ and all pairs of vertices $(a,b)$.

In~\cite{2OrbMani} the authors only use the sets $S_n$ to find the
monodromy $\eta$ of $\M_n$ mapping all the flags of a given facet to
different facets. To ensure that $\eta$ acts this way they only need
to use the fact that $S_n$ is not invariant under non-trivial
automorphisms, and they do not consider any other properties. However,
for our purposes this condition is not enough. We also need
$\eta$ to map the flags of the base facet ``very far away''. This is to
ensure that ``close facets'' have the same base edge as the base
facet. This is why we have proved
Corollary~\ref{c:AristaBase}.

First we will prove that the intersection condition is satisfied for
$k>1$ for every 2-vertex premaniplex $2^{n+1}_I$ where
$0,n\in \ol{I}$.

\begin{teo}\label{t:k>1}
  Let $\X$ be an $(n+1)$-premaniplex with two vertices and with links
  of color 0 and $n$, that is $\X=2^{n+1}_I$ with $0,n\in \ol{I}$. Let
  $\xi$ be the voltage assignment defined in~(\ref{eq:xi}). Then, for
  $k>1$ and for all $m\in [0,n]$ we have that
  \[
    \xi(\Pi^{a,b}_{[0,m]}(\X)) \cap \xi(\Pi^{a,b}_{[k,n]}(\X)) =
    \xi(\Pi^{a,b}_{[k,m]}(\X)),
  \]
  for all pairs of vertices $(a,b)$ in $\X$.
\end{teo}

Before proceeding with the proof, let us introduce a few new
concepts. Let $\M$ be a maniplex of rank $n$ and let $\mu$ be a
flag-permutation of $\M$. Let $i \in \{0,1,\ldots,n-1\}$ and let $F$
and $G$ be $i$-faces of $\M$. If $\mu$ maps every flag with $i$-face
$F$ to a flag with $i$-face $G$ we will say that {\em $\mu$ maps $F$
  to $G$}. In the case where $F=G$ we will say that {\em $\mu$ fixes
  $F$}. If $\mu$ fixes $F$ for every $i$-face $F$ we will say that
{\em $\mu$ fixes $i$-faces}.

To simplify notation, in what follows we will write $\Pi^{a,b}_I$
instead of $\Pi^{a,b}_I(\X)$, taking for granted that we are speaking
about paths on $\X$. Recall that the interval $[k,m]$ is considered to
be the empty set when $k>m$ and that $\xi(\Pi^{a,b}_\emptyset)$ is the
trivial group if $a=b$ and it is the empty set if $a\neq b$.

\begin{dem}
  If $m=n$ there is nothing to prove. If $m<n$ then
  $\xi(\Pi^{a,b}_{[0,m]})$ is generated by monodromies of $\M_n$ mapping
  any white flag $\Phi$ in $\M_n$ to a white flag with the same $i$-faces for $i>m$.

  Now let $W \in \Pi^{a,b}_{[k,n]}$ such that
  $\xi(W) \in \xi(\Pi^{a,b}_{[0,m]})$. Observe that $\xi(W)$ must be a
  monodromy of $\M_n$ that preserves $i$-faces for $i>m$. Let $\Phi$
  be the base flag of the base facet of $\M_n$. Note that, since
  $k>1$, the elements of $\xi(\Pi^{a,b}_{[k,n]})$ are products of flag
  permutations that do not change the 1-face of $\Phi$, that is, for
  every $\omega \in \xi(\Pi^{a,b}_{[k,n]})$ we have that
  $(\Phi\omega)_1 = (\Phi)_1$. This implies that
  $\Phi\omega y_n = \Phi\omega$, as all facets containing $(\Phi)_1$
  must have it as their base edge because of
  Corollary~\ref{c:AristaBaseBase}. Note also that if $1<i<n$, then
  the voltage of all darts of color $i$ is the same. This means that,
  if we write $\xi(W) \in \xi(\Pi^{a,b}_{[k,n]})$ as the product of
  the voltages of the darts of $W$, it acts on $\Phi$ the same way as
  the voltage of the path $W'$ that follows the same colors as $W$ but
  ignoring each occurrence of a dart of color $n$.

  If $W$ uses an even number of darts of color $n$, then $W'$ has the
  same end-points as $W$ and its voltage $\xi(W')$ acts in the same
  way as $\kappa = \xi(W)$ on $\Phi$. Since $\M_n$ is regular, if
  $\kappa$ is a monodromy of $\M_n$ it must coincide with $\xi(W')$,
  since it is also a monodromy acting the same way on some flag. But
  $W'\in \Pi^{a,b}_{[k,n-1]}$, so then we have that
  $\kappa \in \xi(\Pi^{a,b}_{[0,m]}) \cap \xi(\Pi^{a,b}_{[k,n-1]})$.
  Since $\kappa\in \xi(\Pi^{a,b}_{[k,n-1]})$, it must be a monodromy
  that preserves $i$-faces for $i<k$, and since
  $\kappa \in \xi(\Pi^{a,b}_{[0,m]})$, it also preserves $i$-faces for
  $i>m$. Then, by Lemmas~\ref{l:VoltsCerr} and~\ref{l:VoltsAbs},
  $\kappa \in \xi(\Pi^{a,b}_{[k,m]})$.

  Now we want to prove that if $\xi(W)\in \xi(\Pi^{a,b}_{[0,m]})$ then
  $W$ cannot use an odd number of darts of color $n$, thus proving
  that $\xi(\Pi^{a,b}_{[0,m]}) \cap \xi(\Pi^{a,b}_{[k,n]}) =
  \xi(\Pi^{a,b}_{[k,m]})$.

  Let $(e,0)$ be the 1-face of $\Phi$, $(e_0,0)$ be the 1-face of
  $\Phi^1$ and $(e_1,0)$ be the 1-face of $\Phi^{010}$. Here we are
  thinking of $e,e_0$ and $e_1$ as 1-faces of the polytope $\M_{n-1}$,
  which is naturally isomorphic to any facet of $\M_n$ (see
  Figure~\ref{f:Aristas}).

  \begin{figure}
    \begin{center}
      \includegraphics[width=5cm]{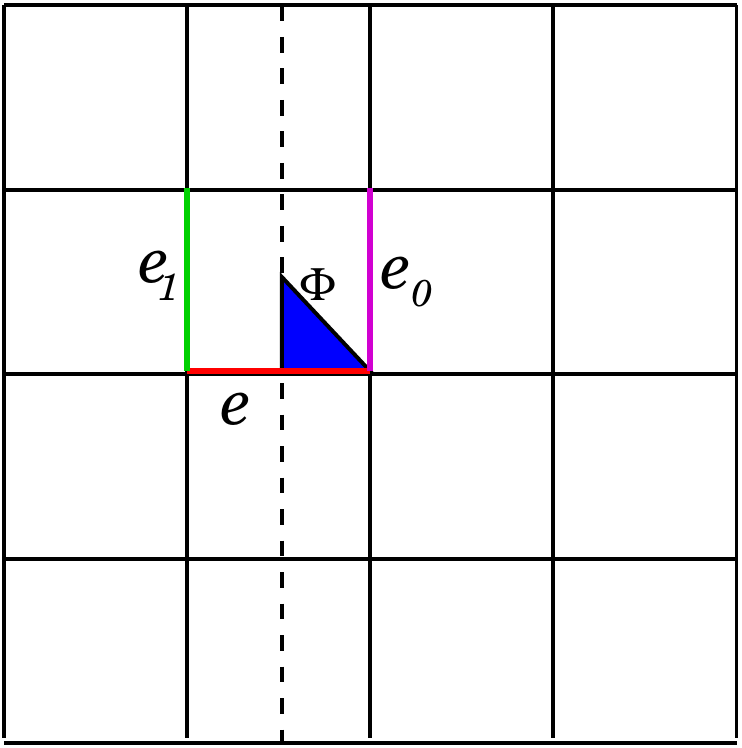}
      \caption{\label{f:Aristas} The edges $e,e_0$ and $e_1$
        illustrated on a 3-face of $\M_{n-1}$.}
    \end{center}
  \end{figure}

  If $F=(f_{n-1},x)$ is a facet of $\M_n$ which has $(e,x)$ as its base
  edge, then $y_n$ interchanges flags with 1-face $(e_0,x)$ with flags
  with 1-face $(e_1,x)$, while it fixes flags with 1-face $(e,x)$ (see
  Figure~\ref{f:ynAristas}).

  \begin{figure}
    \begin{center}
      \includegraphics[width=5cm]{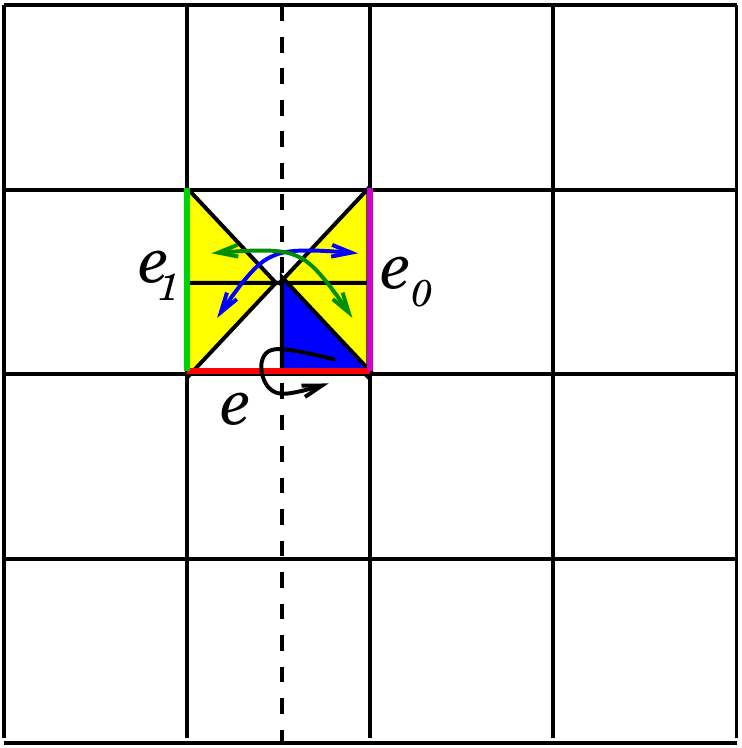}
      \caption{\label{f:ynAristas} If $(e,x)$ is the base edge of a
        facet, then $y_n$ interchanges the edges $(e_0,x)$
        and $(e_1,x)$ while it fixes the edge $(e,x)$.}
    \end{center}
  \end{figure}

  Let $\omega \in \xi(\Pi^{a,b}_{[k,n-1]})$, let $F$ be a facet with
  base edge $(e,x)$, and let $\Psi$ be a flag with facet $F$ and
  1-face $(e_j,x)$ for $j\in\{0,1\}$. Let $(\psi,y) = \Psi\omega$. If
  we write $\omega$ as a product of the voltages of darts, every time
  we change the facet of $\M_n$ (that is, every occurrence of $r_{n-1}$
  or $r_0r_{n-1}$), we must change to a new facet with the same
  edge. This means that the edge of $\Psi\omega$ must be the same as
  the edge of $\Psi$, or in other words, that $(e_j,y)\sim (e_j,x)$
  where $\sim$ is the equivalence relation we used when defining
  $\hat{2}^\p$ for a polytope $\p$. Then $\supp(x+y)$ is contained in
  $\ol{e_j}$.

  Now let $\kappa \in \xi(\Pi^{a,b}_{[k,n]})$ and consider
  $\Phi^1\kappa$ and $\Phi^{10}\kappa$. Since $0\notin I$, then
  $\{\Phi^1\kappa,\Phi^{10}\kappa\}$ has exactly one white flag and
  one black flag. Set $(\psi,x)$ to be the white flag in
  $\{\Phi^1\kappa,\Phi^{10}\kappa\}$. Using again the fact that the
  voltage of a dart with color greater than 1 does not depend on its
  start-point, we may write $\kappa$ as
  $\omega_1 y_n \omega_2 y_n \cdots \omega_{s-1} y_n \omega_s$ where
  $\omega_i \in \xi(\Pi^{a,b}_{[k,n-1]})$ . Each $\omega_i$ may change
  the facet to one where the support of the corresponding vector
  differs in coordinates corresponding to a set contained in
  $\ol{e_0} \cup \ol{e_1}$. Let
  \[
    (\psi_i,x_i) =
    \begin{cases}
      \Phi^1\omega_1 y_n \omega_2 y_n \cdots \omega_i & \text {If }
      (\psi,x)=\Phi^1\kappa,\\
      \Phi^{10}\omega_1 y_n \omega_2 y_n \cdots \omega_i & \text{If }
      (\psi,x)=\Phi^{10}\kappa.
    \end{cases}
  \]
  We claim that $\supp(x_i) \subset \ol{e_0} \cup \ol{e_1}$. This is
  proved by a simple induction on $i$. For $i=1$ we have that $x_1=0$,
  which has support
  $\supp(0)=\emptyset \subset \ol{e_0} \cup \ol{e_1}$. If the claim is
  true for $i$, since $\supp(x_i) \subset \ol{e_0} \cup
  \ol{e_1}$, Corollary~\ref{c:AristaBase} tells us that the facet of $\M_n$
  corresponding to $x_i$ has $(e,x_i)$ as its base edge. Hence, if
  $(\psi_i,x_i)$ has $(e_j,x_i)$ as its 1-face, then
  $(\psi_{i+1},x_{i+1}) = (\psi_i,x_i)y_n\omega_{i+1}$ has $e_{1-j}$
  as its 1-face. This implies that
  $\supp(x_i+x_{i+1}) \subset \ol{e_0} \cup \ol{e_1}$. Then
  \[
    \supp(x_{i+1}) = \supp((x_i+x_{i+1})+x_i) \subset \supp(x_i+x_{i+1})
    \cup \supp(x_i) \subset \ol{e_0} \cup \ol{e_1}.
  \]
  Thus we have proved our claim. Note that $x=x_s$, so
  our claim and Corollary~\ref{c:AristaBase} tell us that the facet
  $(f_{n-1},x)$  of $\M_n$ (where $f_{n-1}$ is the greatest face of
  $\M_{n-1}$) has base edge $(e,x)$.

  Since we assumed that $W$ uses an odd number of darts of color $n$,
  we know by our claim that its voltage $\kappa = \xi(W)$ maps
  $\Phi^1$ or $\Phi^{10}$, with edge $(e_0,0)$, to a flag with edge
  $(e_1,x)$ (for some $x$). But we know also that it maps $\Phi$ to a
  flag with the same edge $(e,0)$. Since $\M_n$ is regular, $\kappa$
  cannot be a monodromy of $\M_n$ (if a monodromy of a regular
  polytope fixes one edge, it must fix all edges), so it cannot be in
  $\xi(\Pi^{a,b}_{[0,m]})$. \qed
\end{dem}

\begin{coro}\label{c:k>1}
  Let $\X$ be an $(n+1)$-premaniplex with two vertices and with links
  of color 0 and $n$, that is $\X=2^{n+1}_I$ with $0,n\in \ol{I}$. Let
  $\xi'$ be the voltage assignment defined in~(\ref{eq:xi'}). Then, for $k>1$ and for all $m\in [0,n]$
  we have that
  \[
    \xi'(\Pi^{a,b}_{[0,m]}) \cap \xi'(\Pi^{a,b}_{[k,n]}) =
    \xi'(\Pi^{a,b}_{[k,m]}),\]
  for all pairs of vertices $(a,b)$ in $\X$.
\end{coro}
\begin{dem}
  Since for $m=n$ there
  is nothing to prove, let us assume that $m<n$.  Let
  $\os{\omega}\in \xi'(\Pi^{a,b}_{[0,m]}) \cap
  \xi'(\Pi^{a,b}_{[k,n]})$. Let $\pi_{\M_n}\omega=\os{\omega}$, that
  is, $\omega$ is the same function as $\os{\omega}$ but considering
  only its action on the $(\M_n)_w$-coordinate. By Theorem~\ref{t:k>1}
  $\omega \in \xi(\Pi^{a,b}_{[k,m]})$. We also know that since $m<n$
  and $\os{\omega}\in \xi'(\Pi^{a,b}_{[0,m]})$ then $\os{\omega}$ is a
  monodromy of $\M_n$ and does not change the $\Z_{2\ell}$-coordinate of
  the elements of $(\M_n)_w\times \Z_{2\ell}$.  Let
  $W\in\Pi^{a,b}_{[k,m]}$ be a path such that $\xi(W)=\omega$. Then
  $\pi_{\M_n}(\xi'(W)) = \xi(W) = \omega = \pi_{\M_n}\os{\omega}$, but
  also $\pi_{\Z_{2\ell}}(\xi'(W)) = Id_{\Z_{2\ell}} =
  \pi_{\Z_{2\ell}}\os{\omega}$. Therefore
  $\xi'(W) = \os{\omega} \in \xi'(\Pi^{a,b}_{[k,m]})$. \qed
\end{dem}

\begin{figure}
  \begin{center}
    \includegraphics[width=5cm]{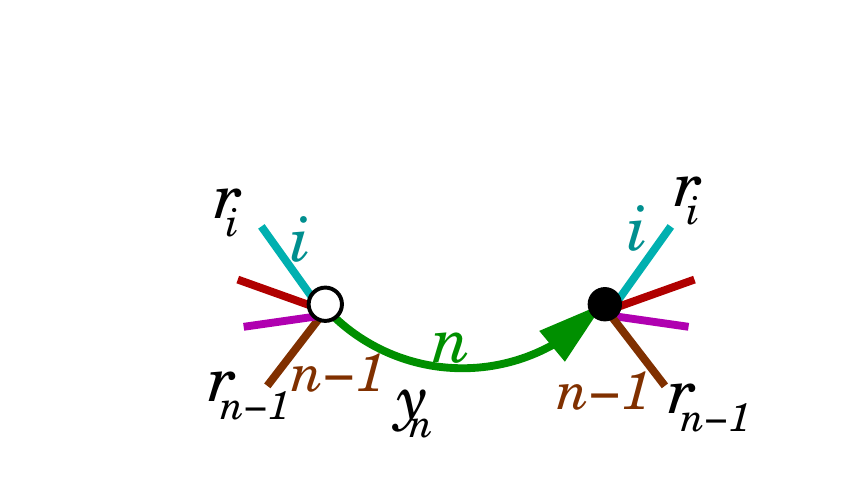}
    \caption{\label{f:Gen2links} The premaniplex
      $(2^{n+1}_{[1,n-1]})_{\ol{0}}$ with its voltage assignment.
    }
  \end{center}
\end{figure}

Now we want to prove that we can remove the restriction $k>1$ from
Theorem~\ref{t:k>1} when $I=[1,n-1]$. First we need the following
lemmas:

\begin{lema}\label{l:Gen2Links}
  Let $a$ be the white vertex on the premaniplex $2^n_{[1,n-1]}$
  (Figure~\ref{f:2links}). Then, the group $\xi(\Pi^a_{[1,n]})$ is
  generated by
  $\{r_i\}_{i=1}^{n-1}\cup \{\os{\rho}_0 r_{n-1} \os{\rho}_0\}$.
\end{lema}
\begin{dem}
  When we remove the link of color 0, we are left with a unique
  spanning tree: the one consisting only of the link of color $n$ and
  voltage $y_n = r_0\os{\rho}_0$ (see Figure~\ref{f:Gen2links}). Then,
  the generator corresponding to a semi-edge of color $i\in [1,n-1]$
  on the white vertex is its voltage $r_i$. The generator
  corresponding to a semi-edge $e$ on the black vertex is the voltage
  of the path $a_nea_n^{-1}$, where $a_n$ is the dart from the white
  vertex to the black one with color $n$. This voltage is
  \[
    \xi(a_nea_n^{-1}) = \xi(a_n)^{-1}\xi(e)\xi(a_n) =
    (y_n)(r_0r_ir_0)(y_n) = \os{\rho}_0r_i\os{\rho}_0.
  \]

  So
  $\xi(\Pi^a_{[1,n]}) =\gen{\{r_i\}_{i=1}^{n-1}\cup \{\os{\rho}_0
    r_{i} \os{\rho}_0\}_{i=1}^{n-1}}$. But since $\os{\rho}_0$ acts as
  an automorphism in each facet, for $i<n-1$ we get
  $\os{\rho}_0r_i\os{\rho}_0 = r_i$, and we get
  $\xi(\Pi^a_{[1,n]}) =\gen{\{r_i\}_{i=1}^{n-1}\cup \{\os{\rho}_0
    r_{n-1} \os{\rho}_0\}}$.\qed
\end{dem}

\begin{lema}\label{l:FijaVertice}
  Let $(e,0)$ be the base edge of the base facet of $\M_n$, and let $u$
  and $v$ be the vertices of $\M_{n-1}$ incident to $e$. If
  $x\in \Z_2^{\Fac(\M_{n-1})}$ is such that
  $\supp(x) \subset \ol{u} \cup \ol{v}$ and $\psi$ is a flag in
  $\M_{n-1}$ with vertex $u$, then
  $(\psi,x)\os{\rho}_0 r_{n-1} \os{\rho}_0 = (\psi, y)$ for a vector
  $y$ satisfying $\supp(y) \subset \ol{u} \cup \ol{v}$.

  Moreover, if $\omega \in \xi(\Pi^a_{[1,n]})$ and
  $(\psi,x)\omega = (\psi',y)$ then $\psi'_0 = u$ and
  $\supp(y) \subset \ol{u} \cup \ol{v}$.
\end{lema}
\begin{dem}
  We know by Corollary~\ref{c:AristaBase} that the facet $(f_{n-1},x)$
  of $\M_n$ has base edge $(e,x)$. Let $\Psi = (\psi,x)$. The vertex
  of $\Psi$ is $(u,x)$. Let $(\psi',x) := \Psi\os{\rho}_0$. We must
  have that $(\psi')_0 = v$. Then
  $\Psi\os{\rho}_0 r_{n-1} = (\psi', x+\chi_F)$ for some
  $F\in \ol{v}$. Since $\supp(x+\chi_F) \subset \ol{u}\cup \ol{v}$ we
  know that the facet $(f_{n-1}, x+\chi_F)$ has base edge
  $(e, x+\chi_{F})$. This implies that
  $\Psi\os{\rho}_0 r_{n-1}\os{\rho}_0 = (\psi', x+\chi_F) \os{\rho}_0
  = (\psi, x+\chi_F)$. To complete the proof use
  Lemma~\ref{l:Gen2Links} and note that each $r_i$ fixes the vertex of
  any flag.
\end{dem}

\begin{teo}\label{t:k=1}
  If $\X=2^{n+1}_{[1,n-1]}$, then $\X$ is the symmetry type graph of a polytope.
\end{teo}

\begin{dem}
  We will prove that the voltage assignment $\xi$ defined
  in~(\ref{eq:xi}) satisfies that
  \begin{equation}
    \label{IntPropk1}
    \xi(\Pi^{a,b}_{[0,m]}) \cap \xi(\Pi^{a,b}_{[k,n]}) =
  \xi(\Pi^{a,b}_{[k,m]})
  \end{equation}
  for all $k,m\in [0,n]$ and all pairs of vertices $(a,b)$. This will
  imply that $\xi'$ (defined in~(\ref{eq:xi'})) also satisfies this
  intersection property (exactly as in the proof of
  Corollary~\ref{c:k>1}). Then, Theorem~\ref{t:IntProp} will imply
  that both $\X^\xi$ and $\X^{\xi'}$ are polytopal, and since we
  already know that $\X^{\xi'}$ has symmetry type graph $\X$ the
  theorem follows.  Since Theorem~\ref{t:k>1} tells us
  that~(\ref{IntPropk1}) holds for $k>1$ and it trivially holds for
  $k=0$, we only need to prove the case when $k=1$.

  Let us first prove the case for closed paths, that is, we will
  prove that
  \[
    \xi(\Pi^a_{[1,n]}) \cap \xi(\Pi^a_{[0,m]}) = \xi(\Pi^a_{[1,m]}).
  \]
  Once again, for $m=n$ there is nothing to prove, so let $m<n$.

We prove first the case when $m<n-1$. Let
  $\omega \in \xi(\Pi^a_{[1,n]}) \cap \xi(\Pi^a_{[0,m]})$. Since
  $\omega \in \xi(\Pi^a_{[0,m]}) = \gen{r_i|i\leq m}$, it must be a
  monodromy of $\M_n$ and since $m<n-1$, it must fix facets. Let
  $\Phi = (\phi,0)$ be the base flag of the base facet of $\M_n$ and
  let $e=(\phi)_1$, $u=(\phi)_0$ and $v=(\phi^0)_0$ (see
  Figure~\ref{f:Vertices}). Lemma~\ref{l:FijaVertice} implies that
  $\Phi\omega = (\psi,x)$ for some flag $\psi$ in $\M_{n-1}$ and some
  vector $x\in\Z_2^{\Fac(\M_{n-1})}$ satisfying $(\psi)_0 = u$ and
  $\supp(x) \subset \ol{u} \cup \ol{v}$. But since $\omega$ must fix
  facets, $x$ is actually 0. This means that
  $(\Phi\omega)_0 = (\Phi)_0$, and since $\M_n$ is regular, this means
  that $\omega$ fixes all 0-faces. Also because of the regularity of
  $\M_n$, the monodromies in $\xi(\Pi^a_{[0,m]}) = \gen{r_i|i\leq m}$
  that fix 0-faces are $\gen{r_i|1\leq i \leq m}$, but this is
  precisely $\xi(\Pi^a_{[1,m]})$.

  Now, if $m=n-1$ we have a little more work to do. We proceed as in
  the previous case, but we cannot ensure that $x=0$. However, we have
  that $\Phi\omega = (\psi,x)$ and the inclusion
  $\supp(x) \subset \ol{u} \cup \ol{v}$ still holds. By writing the
  2-face of $\Phi$ as $(\Phi)_2=(Q,0)$, we know that $Q$ is a
  square. Let $w$ be the opposite vertex of $u$ in $Q$ and $q$ be the
  opposite vertex to $v$, so that the vertices of $Q$ are $uvwq$ in
  cyclical order (see Figure~\ref{f:Vertices}).

  \begin{figure}
    \begin{center}
      \includegraphics[width=5cm]{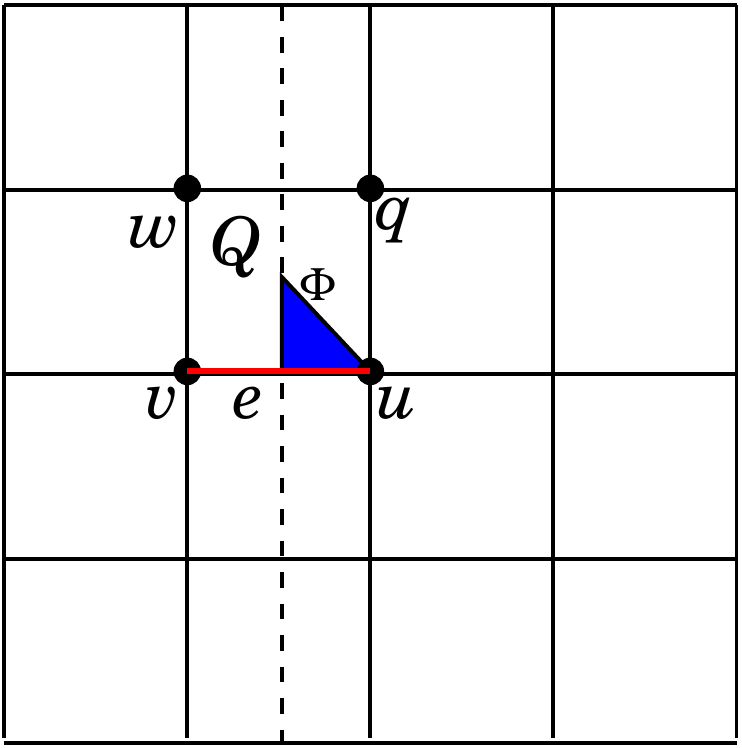}
      \caption{\label{f:Vertices} The 2-face $Q$ and its vertices
        $u,v,w$ and $q$.}
    \end{center}
  \end{figure}

  Since $\omega$ is a monodromy of $\M_n$, it must commute with the
  automorphism $\rho_1$ of $\M_n$, which maps $\Phi$ to $\Phi^1$. This
  implies that
  \[
    \Phi^1\omega = (\Phi\omega)\rho_1 = (\psi, x)\rho_1 = (\psi\rho_1,
    \rho_1 x).
  \]
  Here we have used the way $\Gamma(\M_{n-1})$ acts on $\M_n$
  discussed on Section~\ref{s:2alaM}. Note that we have used
  the same symbol $\rho_1$ to denote an automorphism of $\M_n$ and
  also an automorphism of $\M_{n-1}$, but because of the way
  $\Gamma(\M_{n-1})$ acts on $\M_n$ this is actually not ambiguous.

  Notice that $\supp(\rho_1 x) = (\supp(x)) \rho_1$ must be contained
  in $(\ol{u} \cup \ol{v})\rho_1 = \ol{u} \cup \ol{q}$. On the other
  hand, the vector corresponding to $\Phi^1$ is 0, so
  Lemma~\ref{l:FijaVertice} tells us that
  $\supp(\rho_1 x) \subset \ol{u} \cup \ol{v}$. Every facet incident
  to both $v$ and $q$ must be incident to $Q=v\vee q$ (see
  Corollary~\ref{c:lattice}), and hence also to $u$. Then, the
  intersection $ (\ol{u} \cup \ol{q}) \cap (\ol{u} \cup \ol{v})$ is
  just $\ol{u}$ which means that $\omega$ fixes the vertex of
  $\Phi^1$, and therefore it must fix the vertex of every flag (again,
  $\M_n$ is regular). Hence, we conclude that
  $\xi(\Pi^a_{[1,n]}) \cap \xi(\Pi^a_{[0,m]}) = \xi(\Pi^a_{[1,m]})$.

  Now let us solve the case for open paths, that is, let us prove that
  $\xi(\Pi^{a,b}_{[0,m]}) \cap \xi(\Pi^{a,b}_{[1,n]}) =
  \xi(\Pi^{a,b}_{[1,m]})$ when $a$ is the white vertex and $b$ the
  black vertex.

  First notice that since we only have links of colors 0 and $n$, we
  know that $\xi(\Pi^{a,b}_{[1,m]})$ is in fact empty. So what we
  really want to prove is that there are no monodromies in
  $\xi(\Pi^{a,b}_{[1,n]})$.

  Notice that $\Pi^{a,b}_{[1,n]} = \Pi^a_{[1,n]}a_n$, where once
  again, $a_n$ is the dart of color $n$ from $a$ to $b$. Then
  \[
  \xi(\Pi^{a,b}_{[1,n]}) = \xi(\Pi^a_{[1,n]}a_n) =
  y_n\xi(\Pi^a_{[1,n]}).
  \]

  Let $\omega \in \xi(\Pi^a_{[1,n]})$. We want to prove that
  $y_n\omega$ is not a monodromy of $\M_n$. If it was, it would also
  act as a monodromy on $\M_{n-1}$ (just ignore the
  $\Z_2^{\Fac(\M_{n-1})}$-coordinate). Take the base flag $\phi$ of
  $\M_{n-1}$. Then $\phi y_n\omega = \phi\omega$,
  and Lemma~\ref{l:FijaVertice} tells us that this must be a flag with
  vertex $u$ (the same vertex as $\phi$). Now, if we consider $\phi^1$
  (with vertex $u$) we get that
  $\phi^1y_n\omega = \phi^1\os{\rho}_0 r_0 \omega = \phi^{010}\omega$
  is a flag with vertex $w$. Then $y_n\omega$ cannot act as a
  monodromy on $\M_{n-1}$, since it fixes the vertex of some flags but
  it changes the vertex of others and $\M_{n-1}$ is
  regular. Therefore, $y_n\omega$ is not a monodromy of $\M_n$. We
  have proved that
  $\xi(\Pi^{a,b}_{[0,m]}) \cap \xi(\Pi^{a,b}_{[1,n]}) = \emptyset$.

  Thus we have proved that~(\ref{IntPropk1}) holds. By doing a proof
  analogous to that of Corollary~\ref{c:k>1} we get that
  $\xi'(\Pi^{a,b}_{[0,m]}) \cap \xi'(\Pi^{a,b}_{[k,n]}) =
  \xi'(\Pi^{a,b}_{[k,m]})$ for all $k,m\in
  [0,n]$. Theorem~\ref{t:IntProp} then implies that $\X^{\xi'}$ is
  polytopal. Since we already knew that $\T(\X^{\xi'}) \cong \X$, we
  have found a polytope whose symmetry type graph is isomorphic to $\X$. \qed
\end{dem}

Recall once again that examples of polyhedra (rank 3) of all possible
2-orbit symmetry types exist, in particular those with 1 or 2
links. Applying the constructions $2^\p$ or $\hat{2}^\p$ repeatedly to
the examples of 2 links of Theorem~\ref{t:k=1} or the examples of rank
3 one gets that all symmetry types with 1 or 2 links exist in any rank
higher than two. In addition to
this, as previously discussed, in~\cite{Quirales} Pellicer proves that
chiral polytopes (those with symmetry type $2^n_\emptyset$) exist
in rank 3 and higher, and if we apply the
constructions $2^\p$ or $\hat{2}^\p$ to those, we get all symmetry
types where all the links have consecutive colors. In conclusion, we
have the following theorem:
\begin{teo}\label{t:TiposExisten}
  Let $n\geq 3$ and let $\X=2^n_I$ be a 2-vertex premaniplex of rank
  $n$. Let $\ol{I}:=\{0,1,\ldots,n-1\}\setminus I$ be the set of the
  colors of the links of $\X$. Then, in any of the
  following cases, $\X$ is the symmetry type graph of a polytope.
  \begin{itemize}
  \item $\ol{I}$ has exactly 1 or 2 elements.
  \item $\ol{I}$ is an interval $[k,m] = \{k,k+1,\ldots,m\}$.
  \end{itemize}
\end{teo}

Theorem~\ref{t:TiposExisten} ensures that of the total of $2^n-1$
premaniplexes of rank $n$ with 2 vertices, at least $n^2-n+1$ are
the symmetry type of a polytope ($n$ with 1 link, $\frac{n(n-1)}{2}$
with 2 links and $\frac{n(n-1)}{2}-n+1$ with an interval of links of
size at least 3). It appears that there is still a long way to go.
Nevertheless, Theorem~\ref{t:k>1} and Corollary~\ref{c:k>1} ensure
that to prove that there are polytopes with symmetry type $2^n_I$, one
should only check that ~(\ref{IntPropk1}) is satisfied for $k=1$ for
the voltage $\xi$ (and this would imply that it is also satisfied for
$\xi'$). Sadly, the proof of Theorem~\ref{t:k=1} cannot be easily
generalized to arbitrary $I$. This is because the voltage of a link of
color $1\leq i <n$ would be $r_0r_i$ and this affects most arguments
used, mainly because there would be paths whose voltage change 0-faces
but they do not use the color 0. The author of this paper conjectures
that if this challenge is solved for one example $2^n_I$ with 3 or
more links of non-consecutive colors, the same solution must work for
all others, and that would prove that there exist 2-orbit polytopes of
any possible symmetry type (and rank $n\geq 3$).

\section*{Acknowledgments}
The author would like to thank Gabe Cunningham, Antonio Montero,
Daniel Pellicer, for their helpful comments on early versions of this
paper. The author also gives special thanks to Isabel Hubard, for her
guidance during the writing of the thesis that gave birth to this
paper (and beyond).

The research for this paper was done with the financial support of
CONACyT grant A1-S-21678, and PAPIIT- DGAPA grant IN109023. This paper
was completed while the author had a postdoctoral position at the
Department of Mathematics of Northeastern University.

\bibliographystyle{plain}
\bibliography{Bibliografia}

\end{document}